\documentclass[11pt,a4paper]{article}

\usepackage[utf8]{inputenc}
\usepackage[T1]{fontenc}
\usepackage{amsmath, amssymb}
\usepackage{graphicx}          
\usepackage{subcaption}        
\usepackage{wrapfig}      
\usepackage{booktabs}     
\usepackage{url}         
\usepackage[margin=0.5in]{geometry}
\usepackage[labelfont=bf, singlelinecheck=false]{caption}
\usepackage{hyperref}

\hypersetup{
    colorlinks=true,
    linkcolor=blue,
    filecolor=magenta,      
    urlcolor=cyan,
}

\title{\textbf{Box Allocation Optimization in Meal Kit Delivery}}

\author{
  Thi Minh Thu Nguyen\thanks{Department of Decision Analytics and Risk, Southampton Business School, University of Southampton, Southampton, United Kingdom. Corresponding author. Email: \href{mailto:tmtn1n23@soton.ac.uk}{tmtn1n23@soton.ac.uk}}
  \and
  Loïc Genest\thanks{Data Science Manager, London, United Kingdom.}
  \and
  Alain Zemkoho\thanks{School of Mathematical Sciences, University of Southampton, Southampton, United Kingdom.}
}

\date{September 2025}

\begin{document}

\maketitle
\vspace{-2em}
\begin{abstract}
This study introduces the Box Allocation Problem (BAP), a novel optimization challenge in the \$1.4 billion UK meal kit delivery market. BAP involves assigning orders across multiple production facilities to minimize daily recipe variations while adhering to capacity and eligibility constraints over a 15-day planning horizon. We formulate BAP as a mixed-integer linear programming (MILP) problem and systematically compare the performance of the COIN-OR Branch and Cut (CBC) solver with heuristic methods, including Tabu Search and Iterative Targeted Pairwise Swap. Scalability experiment on instances with up to 100,000 orders show that CBC consistently achieves optimal solutions in under two minutes, maintaining optimality even under dynamic conditions with fluctuating factory capacities and changing customer orders. By reducing day-to-day recipe discrepancies, this approach supports more accurate ingredient forecasting, decreases food waste, and improves operational efficiency across multi-factory network. These results provide the first comprehensive solution framework for temporal allocation problems in meal kit delivery operations.
\end{abstract}


\vspace{1em}
\noindent
\textbf{Abbreviations:} B\&B, Branch and Bound; BAP, Box Allocation Problem; BPP, Bin Packing Problem; CBC, COIN-OR Branch and Cut; CVRP, Capacitated Vehicle Routing Problem; ITPS, Iterative Targeted Pairwise Swap; LD, Lead Day; MILP, Mixed-Integer Linear Programming; TS, Tabu Search; WMAPE, Weighted Mean Absolute Percentage Error.
\vspace{1em}


\section{Introduction}\label{sec1}
The UK meal kit delivery market has experienced rapid expansion, driven by consumer demand for convenient, nutritious, and personalized meal solutions. Reaching \$1.4 billion in revenue by 2024, this sector presents a highly competitive environment where both established food companies and technology-driven startups compete for market share \cite{Statista2024}. Operational excellence in this segment depends heavily on sophisticated supply chain networks capable of coordinating fresh ingredient procurement, inventory management, and time-sensitive distribution across multiple production facilities. The inherent complexity stems from managing perishable goods with short shelf lives, unpredictable demand patterns, and stringent quality requirements -- challenges that transform supply chain optimization from a cost consideration into a strategic imperative for survival \cite{Langham2020, ClearSpider2024}. Companies that fail to achieve operational agility in response to demand volatility and supply disruptions face immediate risks to customer satisfaction, food waste costs, and ultimately, business viability \cite{RMSOmega2024}.

\begin{wrapfigure}{r}{0.25\textwidth}
    \vspace{-1.5em}
    \centering
    \includegraphics[width=3.5cm,height=4.5cm]{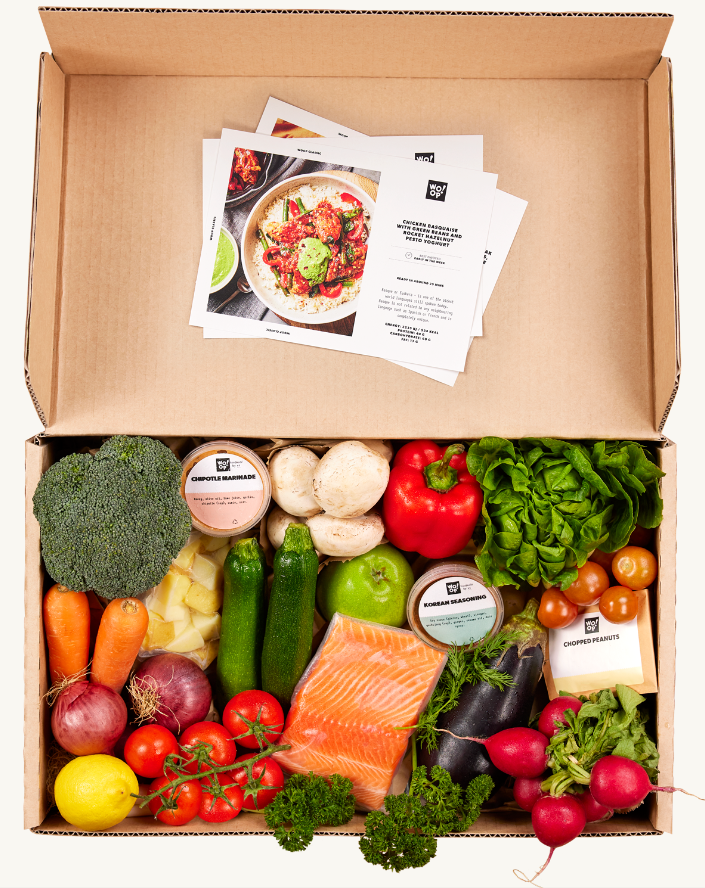}
    \caption{Recipe box~\cite{Woop2024}}
    \label{fig:box}
    \vspace{-1em}
\end{wrapfigure}

Within this operational landscape emerges the Box Allocation Problem (BAP), a critical optimization challenge that involves strategically distributing recipe boxes (Figure \ref{fig:box}) across production facilities to minimize day-to-day variations in recipe allocations. BAP seeks to minimize the Weighted Mean Absolute Percentage Error (WMAPE) of daily recipe distributions while satisfying factory capacity constraint and recipe eligibility restriction \cite{Loic2024}. This objective delivers dual operational benefits that directly impact both sustainability and profitability. First, consistent recipe allocations across facilities enable accurate ingredient demand forecasting, reducing over-procurement and food spoilage -- a critical advantage given that food waste accounts for approximately 10\% of global greenhouse gas emissions \cite{UNFCCC2024}. Second, stabilized allocation patterns ensure reliable ingredient availability across all production sites, minimizing stockouts that can disrupt customer orders and damage brand reputation in this highly competitive market. Despite BAP's fundamental importance to operational efficiency and environmental sustainability in meal kit delivery, it has received minimal attention in the supply chain optimization literature, representing a significant research gap in network-based production planning.

This study addresses the identified research gap by developing a comprehensive mixed-integer linear programming (MILP) formulation of BAP that captures the interplay between facility capacity constraints and recipe eligibility requirements. We systematically benchmark the COIN-OR Branch and Cut (CBC) solver -- an exact optimization method -- against two established heuristic approaches: Iterative Targeted Pairwise Swap (ITPS) and Tabu Search (TS). These heuristics were selected based on their proven effectiveness in structurally similar combinatorial optimization problems, particularly the Bin Packing Problem (BPP) and Capacitated Vehicle Routing Problem (CVRP), which share BAP's fundamental characteristics of resource allocation under capacity constraints \cite{Laporte1987, Lin1965, Gendreau1994}. TS demonstrates particular strength in escaping local optima through its adaptive memory structure and strategic exploration of large solution spaces \cite{Osman1993, Lodi1999}, while CBC guarantees global optimality through systematic branch-and-bound enumeration combined with sophisticated cut generation and primal heuristics \cite{Martello1990, Clausen1999}. Our experimental design encompasses diverse operational scenarios, ranging from 10,000 to 100,000 orders across static and dynamic conditions, including fluctuating factory capacities and evolving customer demand patterns. This comprehensive evaluation framework enables rigorous assessment of each method's computational efficiency, solution quality, and adaptability to real-world complexities.

Our findings highlight the superior performance of the CBC solver in attaining optimal solutions, even under challenging dynamic conditions. By integrating both theoretical foundations and practical applications of BAP, this study makes a substantial contribution to the supply chain optimization literature while offering actionable insights for meal kit delivery enterprises. It bridges a critical knowledge gap by proposing scalable solutions that enhance operational efficiency and promote environmental sustainability in this rapidly growing segment.

The remainder of this paper is structured as follows: Section 2 introduces the characteristics of BAP and presents its mathematical formulation, situating it in relation to two classic optimization problems. Section 3 outlines the optimization methods employed, covering both exact and heuristic approaches. Section 4 reports the results of numerical experiments, comparing the performance of all methods. Finally, Section 5 discusses the implications of the findings, acknowledges limitations, and identifies promising avenues for future research.

\section{Problem Statement}\label{sec2}
\subsection{Problem background}
\begin{figure}[htbp!]
    \centering
    \begin{subfigure}[b]{0.49\textwidth}
        \centering
        \includegraphics[width=\textwidth]{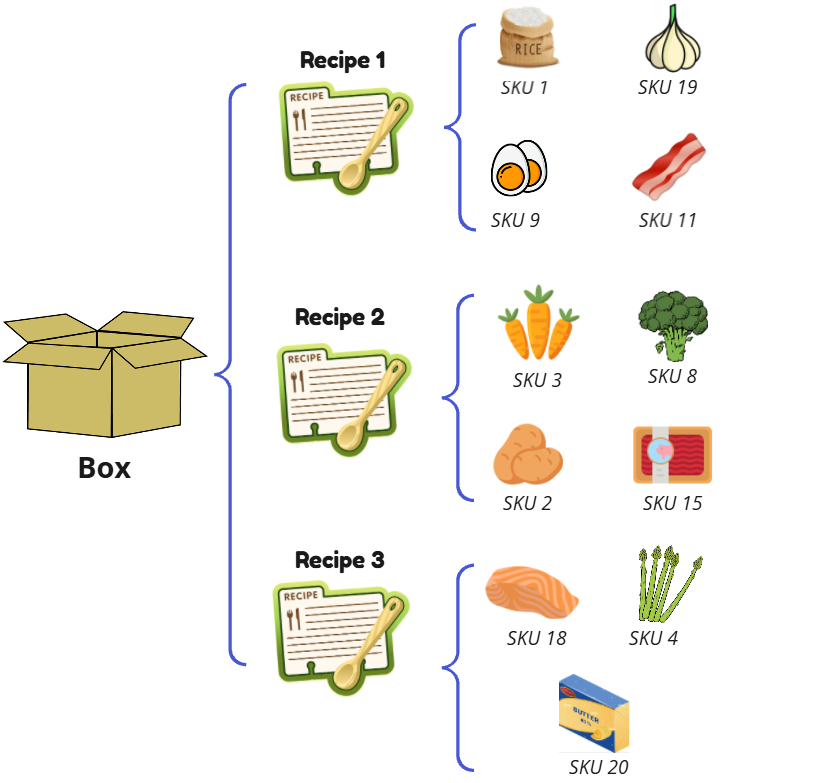}
        \caption{Box structure~\cite{Loic2024}}
        \label{fig:box1}
    \end{subfigure}
    \hfill
    \begin{subfigure}[b]{0.49\textwidth}
        \centering
        \includegraphics[width=\textwidth]{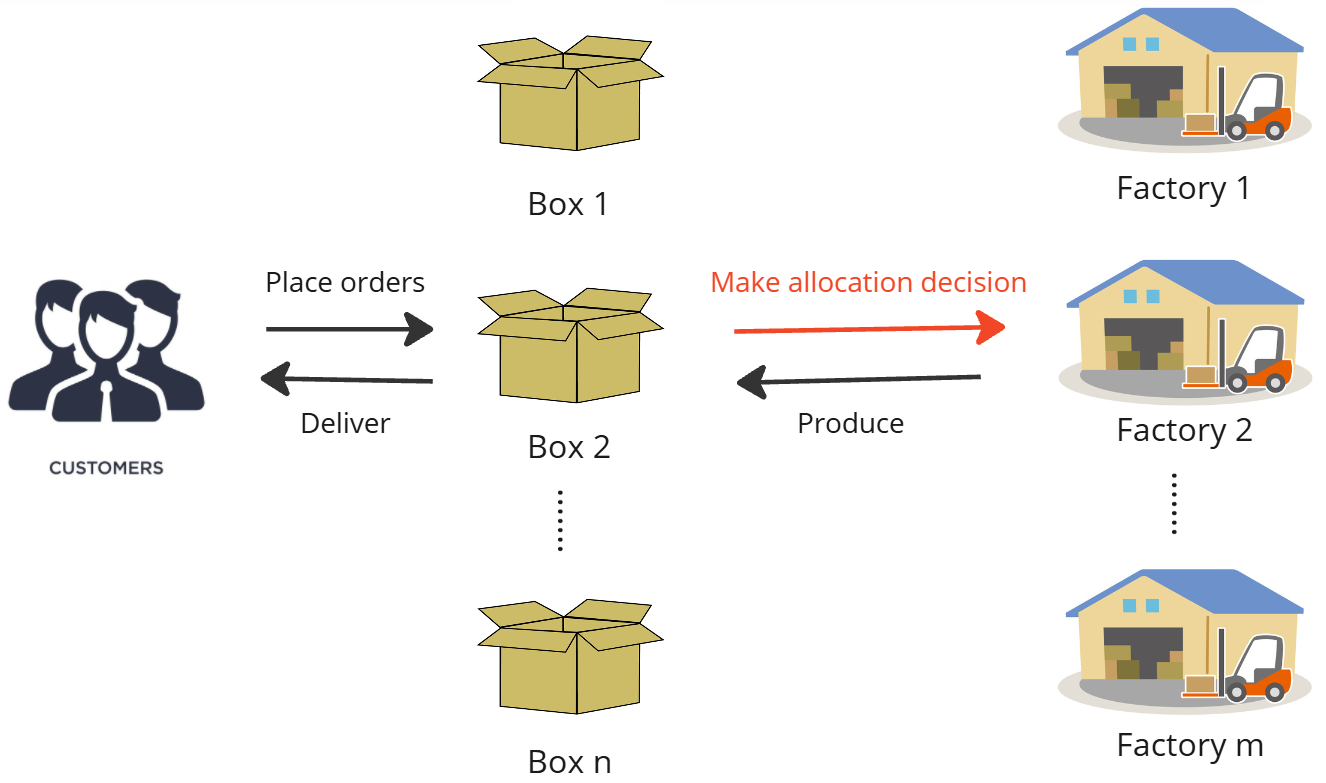}
        \caption{Allocation process~\cite{Loic2024}}
        \label{fig:decision}
    \end{subfigure}
    \caption{Recipe box allocation}
\end{figure}

BAP in meal kit delivery involves assigning recipe boxes, structured as in Figure \ref{fig:box1}, to multiple factories (Figure \ref{fig:decision}). The objective is to minimize daily variations in recipe allocations. For example, if Factory 1 (F1) produced 100 units of Recipe 1 yesterday, the optimal solution would allocate approximately 100 units today.

\begin{figure}[htbp!]
    \centering
    \begin{subfigure}[b]{0.48\textwidth}
        \centering
        \includegraphics[width=\textwidth]{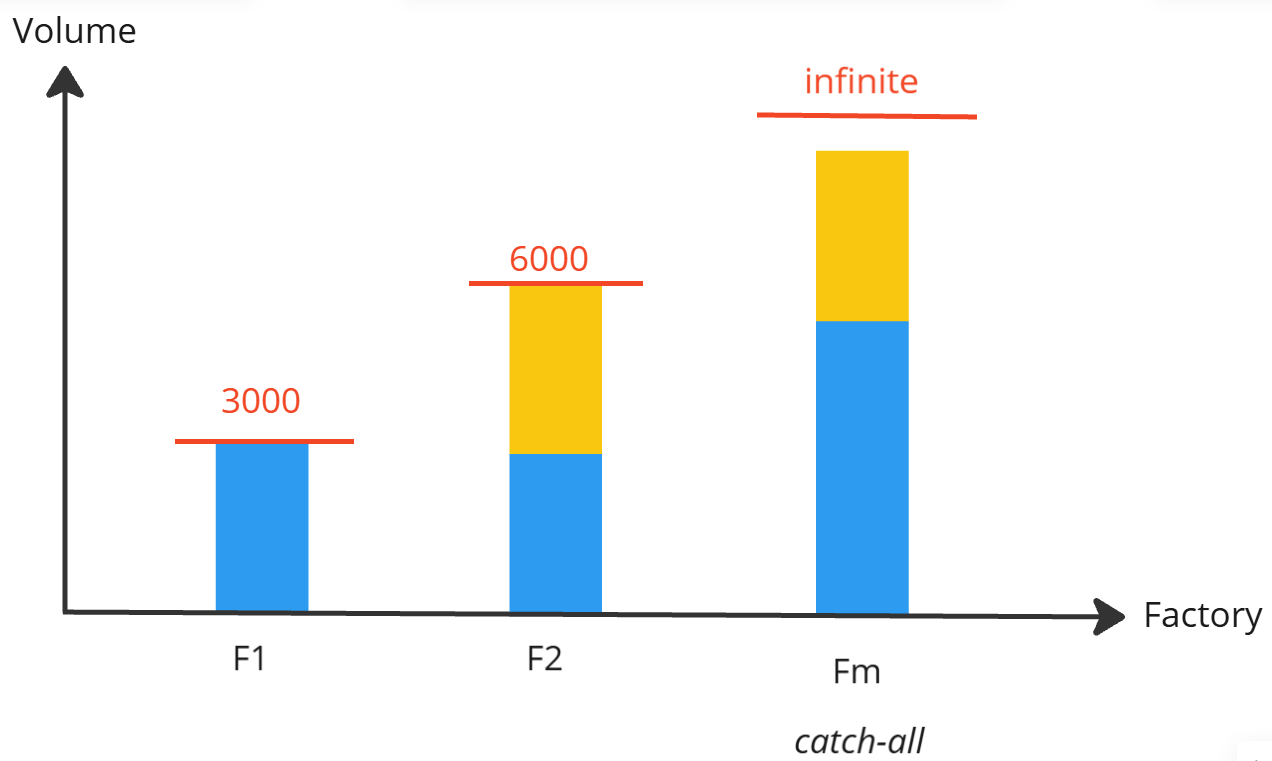}
        \caption{Capacity constraint~\cite{Loic2024}}
        \label{fig:solution}
    \end{subfigure}
    \hfill
    \begin{subfigure}[b]{0.48\textwidth}
        \centering
        \includegraphics[width=\textwidth]{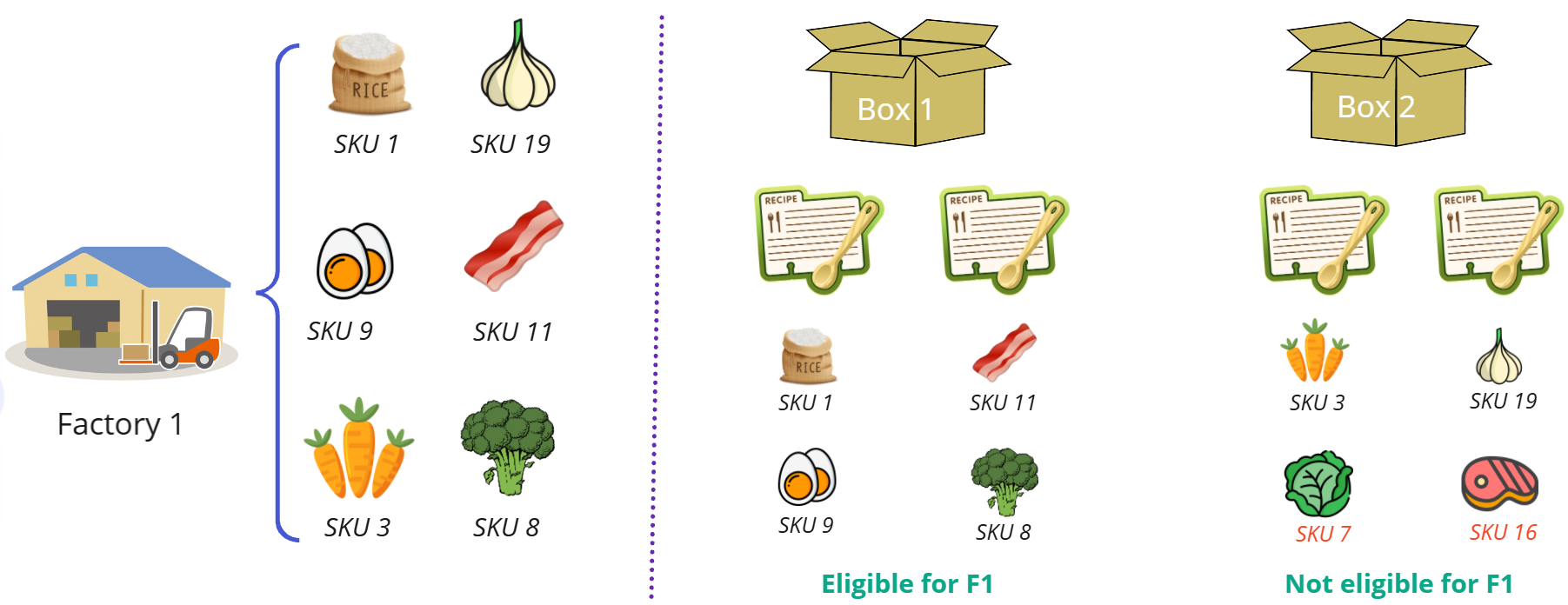}
        \caption{Eligibility constraint~\cite{Loic2024}}
        \label{fig:eligibility}
    \end{subfigure}
    \caption{Main constraints}
\end{figure}

The problem includes two key constraints. First, each factory (except the final simulated factory $m$) has a fixed daily production capacity, limiting the number of boxes it can process. To optimize resource use, solutions must fully meet this capacity (Figure \ref{fig:solution}). Second, not all factories can produce every recipe, so allocation must follow the eligibility constraint in Figure \ref{fig:eligibility}, ensuring recipes are assigned only to capable facilities.

\begin{figure}[htbp!]
    \centering
    \includegraphics[width=0.8\textwidth]{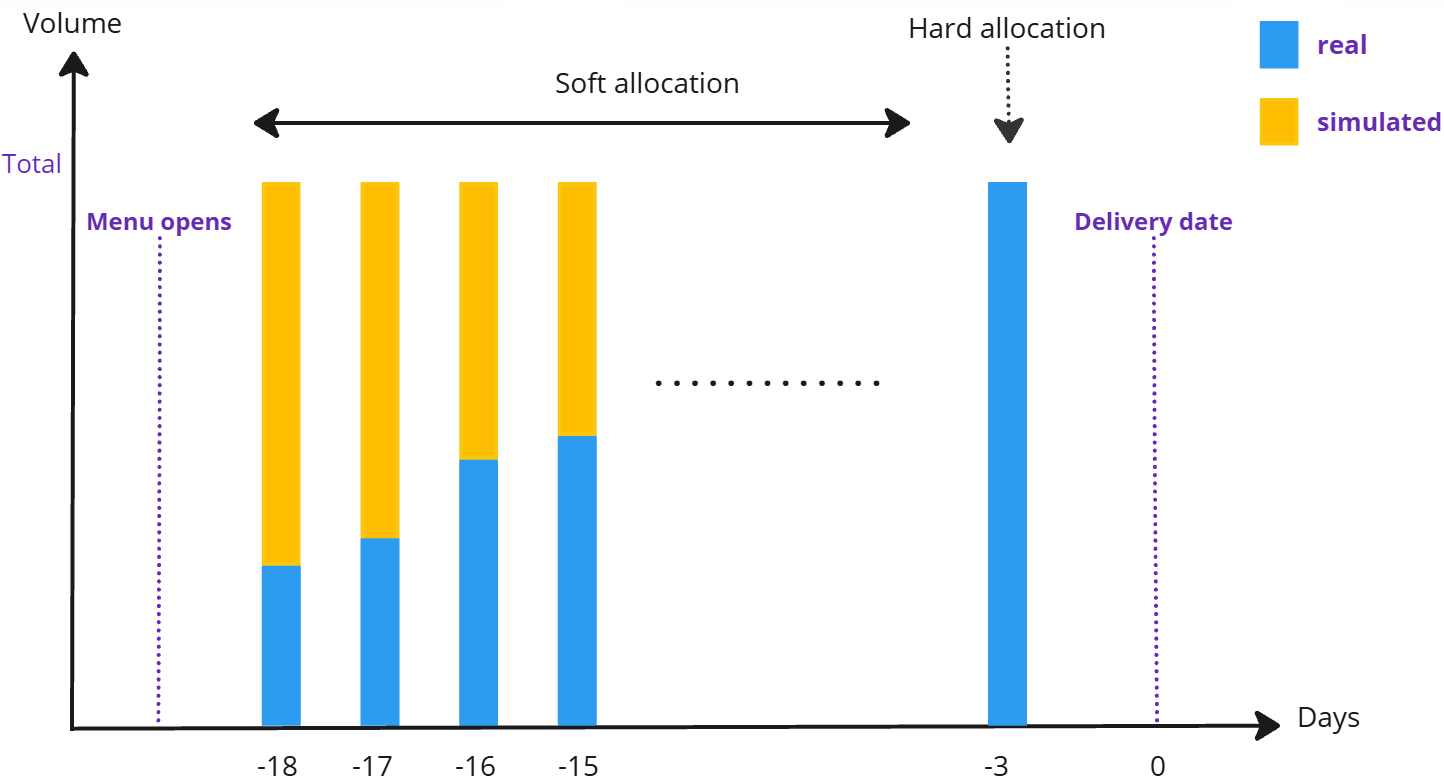}
    \caption{Temporal aspect of BAP~\cite{Loic2024}}
    \label{fig:temporaln}
\end{figure}

A key feature of BAP is its temporal dimension, requiring daily allocation of both actual and simulated orders (Figure \ref{fig:temporaln}). This study uses a 15-day planning horizon from lead day 18 (LD18) to lead day 3 (LD3). The process begins with total demand forecasting, followed by continuous monitoring of actual orders. When real orders fall short of projections, simulated orders are generated to represent expected demand. Initially, both real and simulated orders are distributed across factories (\textit{soft allocation}). As the planning horizon progresses, actual orders replace simulated ones. By LD3, the allocation consists solely of real orders, transitioning to \textit{hard allocation}, which can no longer be changed. The company then has three days to produce and deliver the boxes. This approach ensures stable daily allocations while adapting to fluctuating demand, enabling a smooth shift from forecasts to actual order fulfillment.

\subsection{Mathematical formulation}
WMAPE is a key forecast accuracy metric used in supply chain management and production planning. It measures the sum of absolute differences between actual and forecasted values, divided by the sum of actual values \cite{GeeksforGeeks2021}, with lower values indicating higher accuracy. This research focuses on minimizing WMAPE site in Table~\ref{wmape_comparison} to provide granular insights and assess production consistency at the factory level.

\begin{table}[htbp!]
    \centering
    \caption{Comparison between WMAPE site and WMAPE global}
    \label{wmape_comparison}
    \begin{tabular}{p{0.25\textwidth} p{0.35\textwidth} p{0.35\textwidth}}
        \toprule
        \textbf{Aspect} & \textbf{WMAPE site} & \textbf{WMAPE global} \\ 
        \midrule
        Scope & Considers each factory separately & Aggregates data across all factories \\ 
        Sensitivity & Sensitive to changes in allocation between factories & Not sensitive to changes in allocation between factories \\ 
        Usefulness & Analyzes stability of production plans at each factory & Provides a broader view of recipe distribution changes \\ 
        \bottomrule
    \end{tabular}
\end{table}

The mathematical formulation of BAP employs the following notation:
\begin{itemize}
    \item $ T = \{-18, -17, \dots, -3\} $: Set of allocation days, indexed by $ t $.
    \item $ I = \{1, 2, \dots, n\} $: Set of recipes, indexed by $ i $.
    \item $ J = \{1, 2, \dots, m\} $: Set of factories, indexed by $ j $.
    \item $ O = \{1, 2, \dots, k\} $: Set of orders, indexed by $ o $.
    \item $ a_{i,t-1} $: Quantity of recipe $ i $ at day $ t-1 $.
    \item $ a_{i,t} $: Quantity of recipe $ i $ at day $ t $.
    \item $ a_{i,j,t} $: Quantity of recipe $ i $ allocated to factory $ j $ at day $ t $.
    \item $ a_{i,j,t-1} $: Quantity of recipe $ i $ allocated to factory $ j $ at day $ t-1 $.
    \item $ x_{o,j,t} $: Binary decision variable indicating the allocation of order $ o $ to factory $ j $ on day $ t $.
    \item $ C_{j,t} $: Capacity of factory $ j $ on day $ t $.
    \item $ E_{i,j,t} $: Binary parameter indicating the eligibility of recipe $ i $ at factory $ j $ on day $ t $.
\end{itemize}

We formulate BAP as a MILP problem:
\begin{align}
\underset{x}{\text{Minimize }}&  f(x) = \frac{\sum_{i=1}^{n} \sum_{j=1}^{m} \left| a_{i,j,t} - a_{i,j,t-1} \right|}{\sum_{i=1}^{n} a_{i,j,t}} \label{ProbEq1} \\
 \mbox{subject to}  \;\; &\sum_{o=1}^{k} x_{o,j,t} = C_{j,t}, \qquad \forall j \in J \setminus \{m\},\quad \forall t \in T, \label{ProbEq2} \\
    &x_{o,j,t} \leq \min_{i \in o} E_{i,j,t}, \qquad \forall o \in O, \quad \forall j \in J,\quad \forall t \in T, \label{ProbEq3} \\
    &\sum_{j=1}^{m} x_{o,j,t} = 1, \qquad \forall o \in O, \quad\forall t \in T, \label{ProbEq4} \\
    &x_{o,j,t} \in \{0, 1\},  \qquad \forall o \in O, \quad\forall j \in J, \quad\forall t \in T. \label{ProbEq5}
\end{align}

The objective function \eqref{ProbEq1} minimizes WMAPE site by reducing daily fluctuations in recipe allocations across factories, while accounting for the total recipe volume assigned each day. Recipe quantities $a_{i,j,t}$ are indirectly controlled by the binary decision variable $x_{o,j,t}$ through constraints. When $x_{o,j,t}$ determines which factory j receives order o at time t, it simultaneously establishes how recipes must be distributed. The allocation of an order to a specific factory creates a cascade effect that dictates the required recipe quantities, even though $a_{i,j,t}$ does not explicitly contain $x_{o,j,t}$ in its definition.

Constraint \eqref{ProbEq2} ensures that each factory operates at full capacity to avoid resource waste. Specifically, every factory -- except the simulated catch-all factory $F_m$ -- must receive order quantities that exactly match its production capacity. Constraint \eqref{ProbEq3} enforces the eligibility requirement, allowing an order to be assigned only to a factory capable of producing all recipes required by that order, thereby preventing infeasible manufacturing assignments. Constraint \eqref{ProbEq4} guarantees that each order is assigned to exactly one factory. Finally, Constraint \eqref{ProbEq5} requires the allocation variables to be binary (0 or 1), ensuring that each order is either fully assigned to a particular factory or not assigned at all, reflecting the discrete nature of the assignment decisions.

\vspace{0.5em}
\textbf{WMAPE global} is computed using the following equation:
\begin{equation}
WMAPE_{\text{global}} = \frac{\sum_{i=1}^{n} \left| a_{i,t} - a_{i,t-1} \right|}{\sum_{i=1}^{n} a_{i,t}}, \label{ProbEq6}
\end{equation}
which establishes a theoretical lower bound to evaluate the solution optimality. When decomposed at the site level, this formula becomes:
\begin{equation}
WMAPE_{\text{global}} = \frac{\sum_{i=1}^{n} \left| \sum_{j=1}^{m} a_{i,j,t} - \sum_{j=1}^{m} a_{i,j,t-1} \right|}{\sum_{i=1}^{n} a_{i,t}}, \label{ProbEq7}
\end{equation}
By applying the triangle inequality, we derive the fundamental relationship:
\begin{equation}
WMAPE_{\text{global}} \leq \frac{\sum_{i=1}^{n} \sum_{j=1}^{m} \left| a_{i,j,t} - a_{i,j,t-1} \right|}{\sum_{i=1}^{n} a_{i,t}} (WMAPE_{\text{site}}) \label{ProbEq8}
\end{equation}

Table~\ref{tab:order} presents exemplary order data for two consecutive days, while Table~\ref{tab:errorn} calculates the corresponding WMAPE site and global. This result reveals that despite using the same order data, WMAPE global (0.5) is lower than WMAPE site (1.14), confirming the theoretical relationship in equation \eqref{ProbEq8}.

\begin{table}[htbp!]
    \centering
    \caption{Exemplary orders of LD15 and LD14}
    \label{tab:order}
    \small 
    \setlength{\tabcolsep}{4pt} 
    \begin{tabular}{@{}c@{\hskip 1.2cm}c@{}} 
    \begin{minipage}{0.45\linewidth}
    \centering
    \begin{tabular}{|c|c|c|}
        \hline
        \multicolumn{3}{|c|}{\textbf{LD15}} \\ \hline
        \textbf{Order ID} & \textbf{Recipe IDs} & \textbf{Assigned Factory} \\ \hline
        1 & 1, 10 & F1 \\ \hline
        2 & 2, 3, 5 & F1 \\ \hline
        3 & 4, 6 & F2 \\ \hline
        4 & 2, 3, 7 & F2 \\ \hline
        5 & 3, 5, 9 & F3 \\ \hline
    \end{tabular}
    \end{minipage}
    &
    \begin{minipage}{0.45\linewidth}
    \centering
    \begin{tabular}{|c|c|c|}
        \hline
        \multicolumn{3}{|c|}{\textbf{LD14}} \\ \hline
        \textbf{Order ID} & \textbf{Recipe IDs} & \textbf{Assigned Factory} \\ \hline
        1 & 2, 5 & F1 \\ \hline
        2 & 2, 6, 7 & F3 \\ \hline
        3 & 3, 5, 9 & F1 \\ \hline
        4 & 4, 6, 8 & F2 \\ \hline
        5 & 5, 9, 10 & F2 \\ \hline
    \end{tabular}
    \end{minipage}
    \end{tabular}
\end{table}

\vspace{-1em}
\begin{table}[htbp!]
    \centering
    \caption{WMAPE site and WMAPE global}
    \label{tab:errorn}
    \small 
    \setlength{\tabcolsep}{3pt} 
    \begin{tabular}{@{}c@{\hskip 1.2cm}c@{}} 
    \begin{minipage}{0.48\linewidth}
    \centering
    \begin{tabular}{|c|c|c|c|c|}
        \hline
        \textbf{Recipe} & \textbf{Factory} & \textbf{LD15 ($a_{i,j,t-1}$)} & \textbf{LD14 ($a_{i,j,t}$)} & \textbf{Abs Diff} \\ \hline
        1 & F1 & 1 & 0 & 1 \\ \hline
        2 & F1 & 1 & 1 & 0 \\ \hline
        2 & F2 & 1 & 0 & 1 \\ \hline
        2 & F3 & 0 & 1 & 1 \\ \hline
        3 & F1 & 1 & 1 & 0 \\ \hline
        3 & F2 & 1 & 0 & 1 \\ \hline
        3 & F3 & 1 & 0 & 1 \\ \hline
        4 & F2 & 1 & 1 & 0 \\ \hline
        5 & F1 & 1 & 2 & 1 \\ \hline
        5 & F2 & 0 & 1 & 1 \\ \hline
        6 & F2 & 1 & 1 & 0 \\ \hline
        6 & F3 & 0 & 1 & 1 \\ \hline
        7 & F2 & 1 & 0 & 1 \\ \hline
        7 & F3 & 0 & 1 & 1 \\ \hline
        8 & F2 & 0 & 1 & 1 \\ \hline
        9 & F1 & 0 & 1 & 1 \\ \hline
        9 & F2 & 0 & 1 & 1 \\ \hline
        9 & F3 & 1 & 0 & 1 \\ \hline
        10 & F1 & 1 & 0 & 1 \\ \hline
        10 & F2 & 0 & 1 & 1 \\ \hline
        \multicolumn{3}{|c|}{\textbf{SUM}} & \textbf{14} & \textbf{16} \\ \hline
        \multicolumn{5}{|c|}{\textbf{WMAPE site} = 16/14 = 1.14} \\ \hline
    \end{tabular}
    \end{minipage}
    &
    \begin{minipage}{0.48\linewidth}
    \centering
    \begin{tabular}{|c|c|c|c|}
        \hline
        \textbf{Recipe} & \textbf{LD15 ($a_{i,t-1}$)} & \textbf{LD14 ($a_{i,t}$)} & \textbf{Abs Diff} \\ \hline
        1 & 1 & 0 & 1 \\ \hline
        2 & 2 & 2 & 0 \\ \hline
        3 & 3 & 1 & 2 \\ \hline
        4 & 1 & 1 & 0 \\ \hline
        5 & 2 & 3 & 1 \\ \hline
        6 & 1 & 2 & 1 \\ \hline
        7 & 1 & 1 & 0 \\ \hline
        8 & 0 & 1 & 1 \\ \hline
        9 & 1 & 2 & 1 \\ \hline
        10 & 1 & 1 & 0 \\ \hline
        \multicolumn{2}{|c|}{\textbf{SUM}} & \textbf{14} & \textbf{7} \\ \hline
        \multicolumn{4}{|c|}{\textbf{WMAPE global} = 7/14 = 0.5} \\ \hline
    \end{tabular}
    \end{minipage}
    \end{tabular}
\end{table}

\begin{figure}[htbp!]
    \centering
    \includegraphics[width=0.6\textwidth]{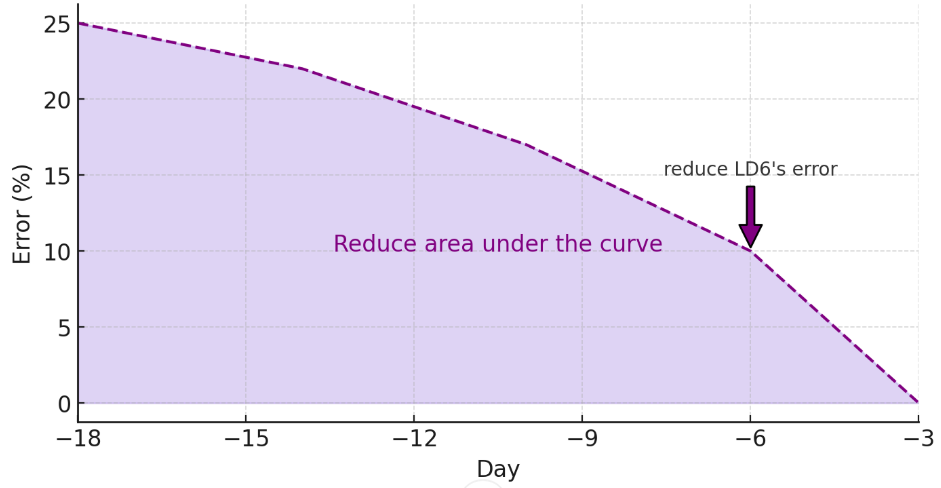}
    \caption{Objective of proxy optimization}
    \label{fig:error}
\end{figure}

The company's ultimate goal is to minimize the area under the curve in Figure \ref{fig:error}, which represents the WMAPE site when comparing the final day's allocation with each of the previous days. However, directly optimizing this metric is infeasible due to forecast uncertainties before final day LD3. These uncertainties arise from uncontrollable factors, including the variability of simulated orders, customer modifications to existing orders, and the continuous arrival of new actual orders throughout the planning period. To tackle this challenge, the company adopts a more practical approach by minimizing allocation variations between consecutive days, as defined in Equation \ref{ProbEq1}. This strategy reframes BAP as a \textit{proxy optimization problem}, where solving a simpler, related problem serves as a substitute for optimizing the more complex objective \cite{Zangl2006}. By reducing day-to-day allocation discrepancies, this approach effectively lowers cumulative error over the entire planning period, making it a feasible way to achieve the company's broader strategic goal.

\begin{figure}[htbp!]
    \centering
    \includegraphics[width=0.7\textwidth]{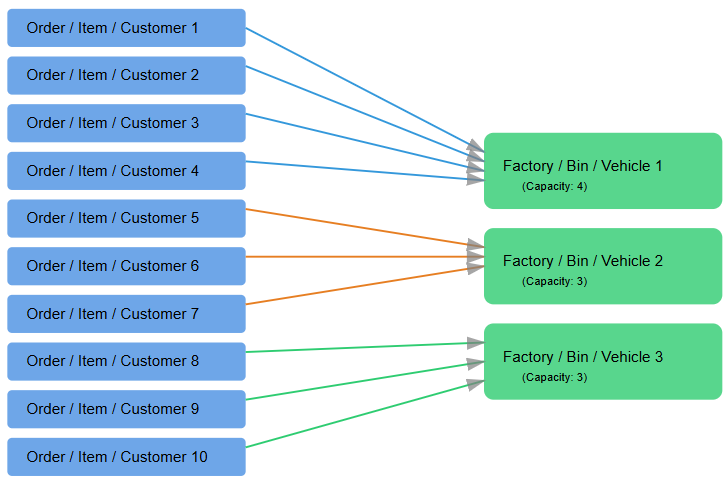}
    \caption{Similarities between BAP, BPP, and CVRP}
    \label{fig:network}
\end{figure}

BAP shares similarities with two well-known combinatorial optimization problems: BPP and CVRP. These problems rely on a network structure, where nodes represent key entities (factories, storage bins, or vehicles) and edges define decision variables (order allocations, item assignments, or routing sequences). Allocations are subject to capacity constraints to ensure efficient resource utilization, as illustrated in Figure \ref{fig:network}. Despite these shared features, each problem has a distinct mathematical formulation. BPP minimizes the number of bins while keeping total item weight or volume within capacity limits. CVRP focuses on reducing total travel costs or distances while ensuring vehicle capacity is not exceeded. BAP introduces the unique objective of minimising temporal variations in recipe allocations across factories.

Their similarities highlight opportunities for cross-domain optimization techniques. The B\&B method from BPP has proven effective across various allocation problems \cite{Martello1990}. Likewise, heuristic approaches such as 2-Opt and Tabu Search from CVRP have influenced allocation strategies in multiple domains \cite{Lin1965, Glover1986}. However, applying these methods to BAP requires substantial modifications to account for its unique constraints and the emphasis on maintaining stable daily allocations. These distinctions place BAP beyond traditional BPP and CVRP models, requiring problem-specific mathematical formulations while strategically leveraging established optimization techniques.

\section{Problem Solving}\label{sec3}
\begin{figure}[htbp!]
    \centering
    \includegraphics[width=\textwidth]{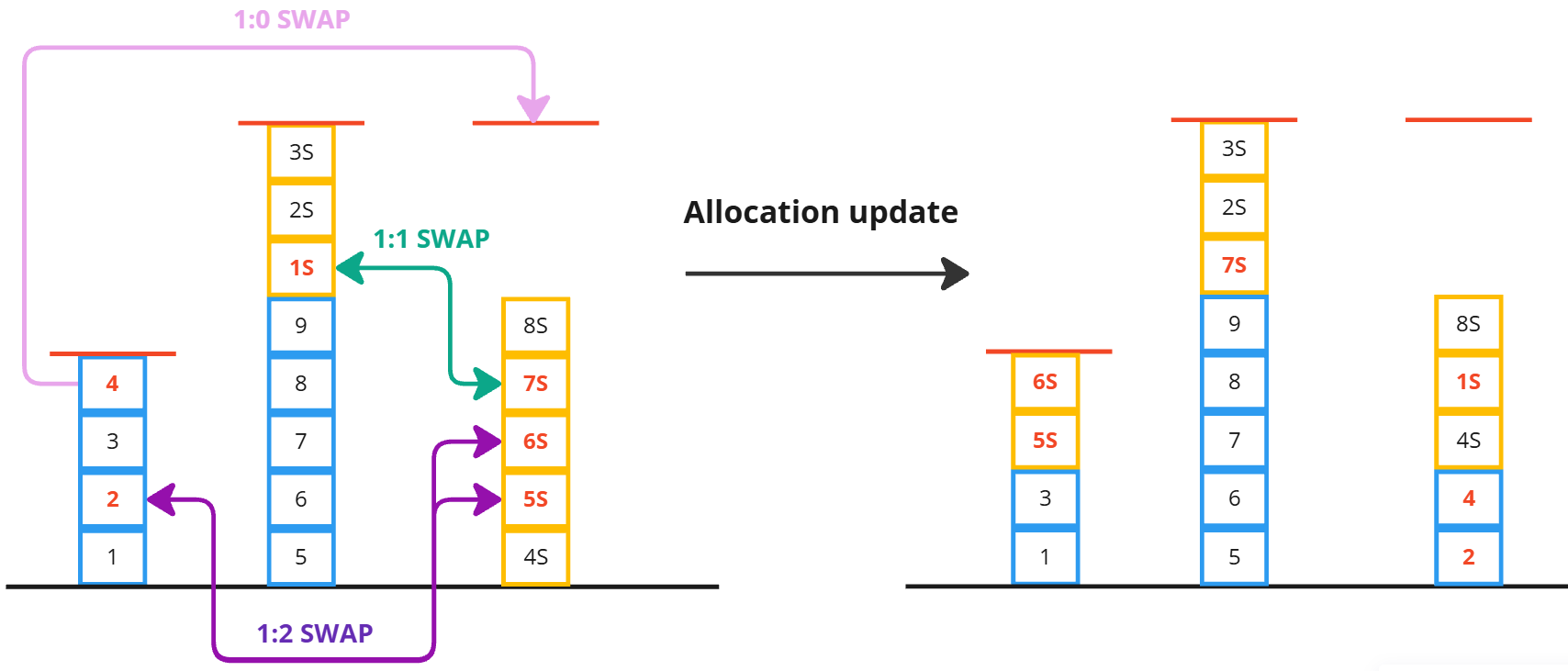}
    \caption{Types of box swaps}
    \label{fig:swap}
\end{figure}

The company seeks an efficient approach to minimize recipe discrepancies between two allocation decisions. Figure \ref{fig:swap} illustrates the different swap types employed in heuristics, where blue boxes represent real orders and orange boxes denote simulated orders. To simplify the process, this study adopts the 1:1 swap method. The solving procedure consists of the following steps:

\begin{itemize}
    \item \textbf{Exact method:} CBC solver is directly applied to current orders to determine the optimal allocation that minimizes recipe discrepancies relative to the previous day's allocation.
    \item \textbf{Heuristics:} The process starts with a greedy algorithm to construct an initial feasible solution. Orders are first allocated to F1, prioritizing those with fewer eligible factories (F1-F3 first, then F1-F2-F3) until capacity is reached. Remaining orders with eligible factories F1-F2-F3 and F2-F3 are allocated to F2 until its capacity is full. Any unallocated orders are assigned to F3, the catch-all factory without constraints. ITPS and TS are then applied to refine the initial allocation.
    \item \textbf{Evaluation:} The performance and scalability of each method are assessed by comparing the final WMAPE site against the WMAPE global and analyzing optimization times for varying order quantities. The most effective method is subsequently applied in temporal tests to evaluate performance over a 15-day planning period and its adaptability to dynamic changes.
\end{itemize}

This article first examines the application of the exact method to address BAP. The CBC solver, developed under the Computational Infrastructure for Operations Research (COIN-OR) project, offers a robust open-source solution for MILP problems \cite{Forrest2005}. Its seamless integration with the COIN-OR linear programming solver and cut-generation library makes it particularly effective for tackling complex optimization challenges. In this study, the CBC solver can solve BAP in just a few steps. For example, consider a hypothetical order quantity of 100,000; the process proceeds as follows:

\begin{enumerate}
    \item The model is loaded and initialized with 210,302 rows, 300,600 columns, and 1,359,124 elements. Preprocessing fixes 110,000 variables, reducing the problem size. Solving the continuous relaxation yields an objective value of 4,646.
    \item The Feasibility Pump heuristic (Table \ref{tab:heu}) is employed to quickly obtain an initial integer-feasible solution. Coefficient Diving (Table \ref{tab:heu}) is then applied, successfully identifying the optimal solution. For most order quantities below 100,000, the Feasibility Pump alone suffices to find the optimal solution immediately.
    \item A limited B\&B process resolves any remaining general integer variables not fixed by heuristics. This step focuses on the reduced problem size, branching only on the unfixed variables.
    \item A mini B\&B is conducted to further refine the solution. While this step does not improve the objective, it serves to validate the quality of the heuristic-based solution.
    \item The process concludes without extensive branching, as the optimal solution is proven at the root node following the limited B\&B and mini B\&B steps.
\end{enumerate}

\begin{table}[htbp!]
    \centering
    \caption{Primal heuristics in CBC solver}
    \label{tab:heu}
    \begin{tabular}{|p{0.15\textwidth}|p{0.4\textwidth}|p{0.4\textwidth}|}
    \hline
    \textbf{Category} & \textbf{Feasibility Pump} & \textbf{Coefficient Diving} \\ \hline
    Overview & Introduced by Fischetti et al. \cite{Fischetti2005}, this heuristic rapidly finds feasible solutions by alternating between linear programming (LP) relaxation and integer solutions. & Based on works by Berthold \cite{Berthold2006}, and Paulus and Krause \cite{Paulus2023}, this heuristic prioritizes variables with minimal constraint violations. \\ \hline
    Steps & 
    \begin{enumerate}
        \item Start with the LP relaxation solution.
        \item Round to the nearest integer-feasible point.
        \item Find the closest LP solution to this integer point.
        \item Repeat until a feasible solution is found or a stopping criterion is met.
    \end{enumerate} &  
    \begin{enumerate}
        \item Select variables based on \textit{locks} (potential constraint violations).
        \item Prioritize variables with minimal locks and smallest fractionality.
        \item Iteratively bound selected variables to guide the solution toward feasibility.
    \end{enumerate} \\ \hline
    Key advantage & The \textit{pumping} action between LP-feasible and integer-feasible solutions enables fast convergence to high-quality feasible solutions. & This targeted approach efficiently finds optimal or near-optimal solutions early in the process. \\ \hline
    \end{tabular}
\end{table}

Regarding the two heuristics, ITPS is inspired by the 2-Opt method commonly applied in VRP \cite{Toth2002}. While 2-Opt swaps adjacent elements in a route to minimize travel distance or cost, ITPS swaps two orders between factories to minimize site-level WMAPE. ITPS has three defining features: it accepts only improving swaps, selects candidates based on recipe differences, and preserves feasibility by swapping only two orders at a time.

TS is a metaheuristic grounded in adaptive memory techniques, employing a memory structure known as the \textit{tabu list} to guide the search more effectively \cite{Glover1997}. The term \textit{tabu} originates from Polynesian cultures, where it denotes prohibitions \cite{Radcliffe1939}. In TS, tabu moves are those temporarily forbidden to prevent cycling and encourage exploration of new regions in the solution space \cite{Glover1989}. A central concept is the \textit{tenure}, which defines the duration a move remains in the tabu list \cite{Gendreau2003}. This mechanism prevents short-term cycling while allowing previously explored moves to be revisited after adequate diversification \cite{Battiti1994}. By balancing intensification (selecting the best move) with diversification (introducing random swaps), TS ensures a comprehensive search of the solution space.

\section{Numerical experiments}\label{sec4}
\subsection{Experimental design}
Simulated data replicating the characteristics of real data is generated to evaluate the effectiveness of three proposed methods. The resulting order list follows realistic rules: 30\% of the total orders are eligible for F1, 60\% for F2, and all orders are eligible for F3, which has no constraints. When the proportion of real orders increases over time, all real orders from the previous day are carried over to the current day with their IDs unchanged, facilitating allocation tracking. Simulated orders are generated anew each day to make up for any shortfall in the total quantity, resulting in dynamic IDs. Each order consists of 1 to 4 recipes chosen from a total of 100 available recipes, with eligibility based on their grouping. Recipes in Group 1 (1–29) are eligible exclusively for F1, Group 2 (30–49) qualifies for both F1 and F2, Group 3 (50–89) is eligible only for F2, and Group 4 (90–100) is eligible solely for F3. An order qualifies for a factory only if all its recipes meet the factory's eligibility. However, for F3-only orders, any recipe can be included as long as the order contains at least one recipe from Group 4. Table~\ref{tab:order1} presents exemplary order data along with solutions that satisfy the capacity constraints: 25\% of the total quantity for F1 and 50\% for F2. Importantly, the allocated factories for four existing real orders on LD12 remain unchanged, ensuring a minimum WMAPE site.

\begin{table}[htbp!]
    \centering
    \caption{Exemplary orders and solutions}
    \label{tab:order1}
    \footnotesize 
    \begin{tabular}{cc}
        \begin{tabular}{|c|c|c|c|}
            \hline
            \multicolumn{4}{|c|}{\textbf{LD12 (46\% real orders)}} \\ \hline
            \textbf{Order ID} & \textbf{Recipe IDs} & \textbf{Is real} & \textbf{Eligible factories} \\ \hline
            1 & 30 & True & F1, F2, F3 \\ \hline
            2 & 8, 5, 24 & True & F1, F3 \\ \hline
            3 & 22 & True & F1, F3 \\ \hline
            4 & 87 & True & F2, F3 \\ \hline
            5 & 52, 51, 55, 63 & False & F2, F3 \\ \hline
            6 & 82, 88 & False & F2, F3 \\ \hline
            7 & 85 & False & F2, F3 \\ \hline
            8 & 84, 76 & False & F2, F3 \\ \hline
            9 & 93, 1, 36, 76 & False & F3 \\ \hline
            10 & 28, 95, 20 & False & F3 \\ \hline
        \end{tabular}
        &
        \begin{tabular}{|c|c|c|c|}
            \hline
            \multicolumn{4}{|c|}{\textbf{LD11 (52\% real orders)}} \\ \hline
            \textbf{Order ID} & \textbf{Recipe IDs} & \textbf{Is real} & \textbf{Eligible factories} \\ \hline
            1 & 30 & True & F1, F2, F3 \\ \hline
            2 & 8, 5, 24 & True & F1, F3 \\ \hline
            3 & 22 & True & F1, F3 \\ \hline
            4 & 87 & True & F2, F3 \\ \hline
            5 & 74 & True & F2, F3 \\ \hline
            6 & 85 & False & F2, F3 \\ \hline
            7 & 89, 73, 86 & False & F2, F3 \\ \hline
            8 & 54, 52 & False & F2, F3 \\ \hline
            9 & 100, 99 & False & F3 \\ \hline
            10 & 91, 13 & False & F3 \\ \hline
        \end{tabular}
    \end{tabular}

    \vspace{1em}
    
    \begin{tabular}{cc}
        \begin{tabular}{|c|c|}
            \hline
            \multicolumn{2}{|c|}{\textbf{Solution for LD12}} \\ \hline
            \textbf{Factory} & \textbf{Allocated orders} \\ \hline
            F1 & 2, 3 \\ \hline
            F2 & 1, 4, 5, 6, 8 \\ \hline
            F3 & 7, 9, 10 \\ \hline
        \end{tabular}
        &
        \begin{tabular}{|c|c|}
            \hline
            \multicolumn{2}{|c|}{\textbf{Solution for LD11}} \\ \hline
            \textbf{Factory} & \textbf{Allocated orders} \\ \hline
            F1 & 2, 3 \\ \hline
            F2 & 1, 4, 5, 7, 8 \\ \hline
            F3 & 6, 9, 10 \\ \hline
        \end{tabular}
    \end{tabular}
\end{table}

To evaluate the performance of the proposed methods, three key indicators are used. First, the optimality gap measures the difference between the MAPE site and global, with WMAPE global serving as the lower bound. CBC solver is expected to match WMAPE site to WMAPE global, as the exact method aims to find the optimal solution. Second, the improvement percentage measures the reduction in WMAPE site relative to the initial solution for heuristics, enabling performance comparison since these methods may not always reach the optimal result. Finally, optimization time is assessed by the duration each method takes to find the best solution, with a maximum limit of 10 minutes (600 seconds).

\subsection{Results}
\subsubsection{Benchmark test}
The following results are based on 10,000 orders over two consecutive days: LD12 (46\% real orders) and LD11 (52\% real orders). The CBC solver uses the LD12 allocation decision to directly determine the best allocation for LD11, minimizing the site-level WMAPE between the two days without swaps. For the heuristics, the number of iterations strongly influences their ability to reach high-quality solutions within reasonable time limits. To identify suitable iteration counts, tests were conducted on 30 different sets of 10,000 orders, resulting in 1,500 iterations for ITPS and 500 iterations for TS as the optimal settings. Using these values, the allocations for LD11 before and after applying ITPS and TS are summarized in Table~\ref{tab:benchmark_results}. The results show that B\&B or CBC solver remains the most effective method, achieving the optimal solution in the shortest time. Between the two heuristics, TS slightly outperforms ITPS in both solution quality and computational speed.

\begin{table}[htbp!]
\centering
\caption{Benchmark test results}
\label{tab:benchmark_results}
\footnotesize
\setlength{\tabcolsep}{2pt} 
\begin{tabular}{p{0.4\textwidth}@{\hspace{0.3em}}p{0.47\textwidth}}
\begin{tabular}{|l|c|}
\hline
\multicolumn{2}{|c|}{\textbf{B\&B}} \\
\hline
\multicolumn{2}{|c|}{\textbf{LD12}} \\
\hline
F1 & 2500 orders, 2500 real \\
\hline
F2 & 5000 orders, 1890 real \\
\hline
F3 & 2500 orders, 210 real \\
\hline
\multicolumn{2}{|c|}{\textbf{LD11}} \\
\hline
F1 & 2500 orders, 2500 real \\
\hline
F2 & 5000 orders, 2593 real \\
\hline
F3 & 2500 orders, 107 real \\
\hline
\textbf{WMAPE site} & 0.054 \\
\hline
\textbf{WMAPE global} & 0.054 \\
\hline
\textbf{Optimization time} & 2.89 seconds \\
\hline
\end{tabular}
&
\begin{tabular}{|l|c|c|}
\hline
\multicolumn{1}{|c|}{} & \multicolumn{1}{c|}{\textbf{ITPS}} & \multicolumn{1}{c|}{\textbf{TS}} \\
\hline
\multicolumn{3}{|c|}{\textbf{LD12}} \\
\hline
F1 & \multicolumn{2}{c|}{2500 orders, 2500 real} \\
\hline
F2 & \multicolumn{2}{c|}{5000 orders, 1890 real} \\
\hline
F3 & \multicolumn{2}{c|}{2500 orders, 210 real} \\
\hline
\multicolumn{3}{|c|}{\textbf{LD11 (Before)}} \\
\hline
F1 & \multicolumn{2}{c|}{2500 orders, 2500 real} \\
\hline
F2 & \multicolumn{2}{c|}{5000 orders, 2422 real} \\
\hline
F3 & \multicolumn{2}{c|}{2500 orders, 278 real} \\
\hline
\multicolumn{3}{|c|}{\textbf{LD11 (After)}} \\
\hline
F1 & 2500 orders, 2500 real & 2500 orders, 2500 real \\
\hline
F2 & 5000 orders, 2436 real & 5000 orders, 2429 real \\
\hline
F3 & 2500 orders, 264 real & 2500 orders, 271 real \\
\hline
\multicolumn{3}{|c|}{\textbf{WMAPE site}} \\
\hline
Before & 0.074 & 0.074 \\
\hline
After & 0.054 & 0.054 \\
\hline
Improvement & 26.21\% & 26.97\% \\
\hline
\multicolumn{3}{|c|}{\textbf{WMAPE global}} \\
\hline
\multicolumn{3}{|c|}{0.054} \\
\hline
\textbf{Optimization time} & 19.45 seconds & 15.35 seconds \\
\hline
\end{tabular}
\end{tabular}
\end{table}

\subsubsection{Scalability test}
\begin{figure}[htbp!]
    \centering
    \begin{minipage}[b]{0.49\textwidth}
        \centering
        \includegraphics[width=\textwidth]{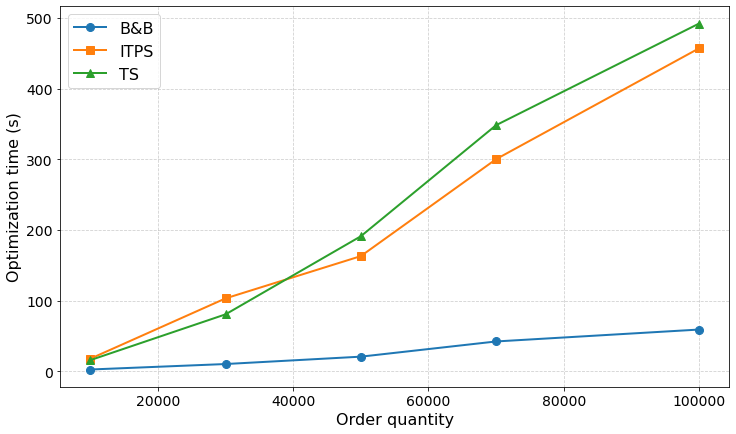}
        \caption{Comparison of optimization time}
        \label{fig:sca}
    \end{minipage}
    \hfill
    \begin{minipage}[b]{0.49\textwidth}
        \centering
        \includegraphics[width=\textwidth]{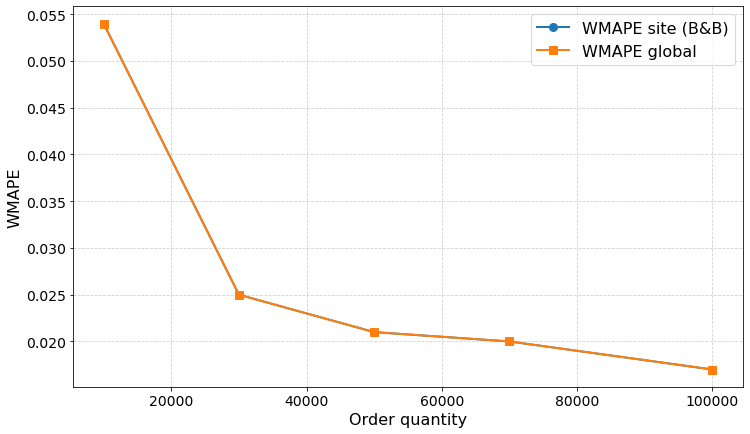}
        \caption{WMAPEs of B\&B}
        \label{fig:sca01}
    \end{minipage}
    
    \vspace{1em}
    \begin{minipage}[b]{0.49\textwidth}
        \centering
        \includegraphics[width=\textwidth]{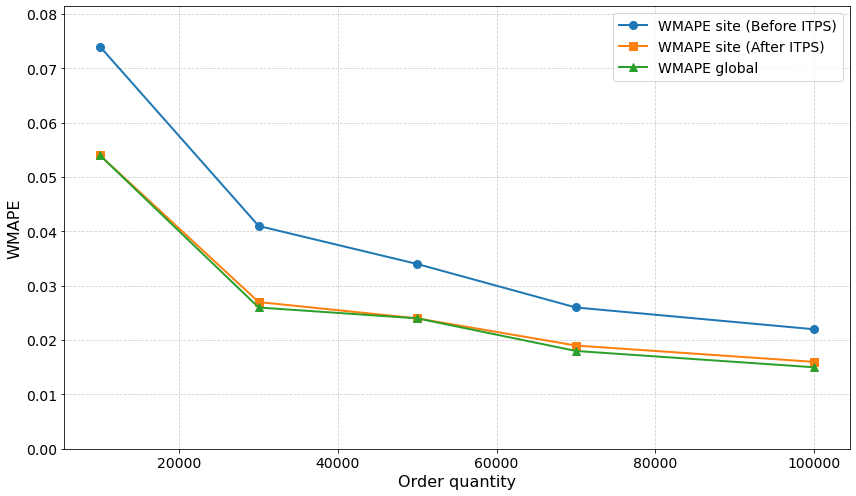}
        \caption{WMAPEs of ITPS}
        \label{fig:sca11}
    \end{minipage}
    \hfill
    \begin{minipage}[b]{0.49\textwidth}
        \centering
        \includegraphics[width=\textwidth]{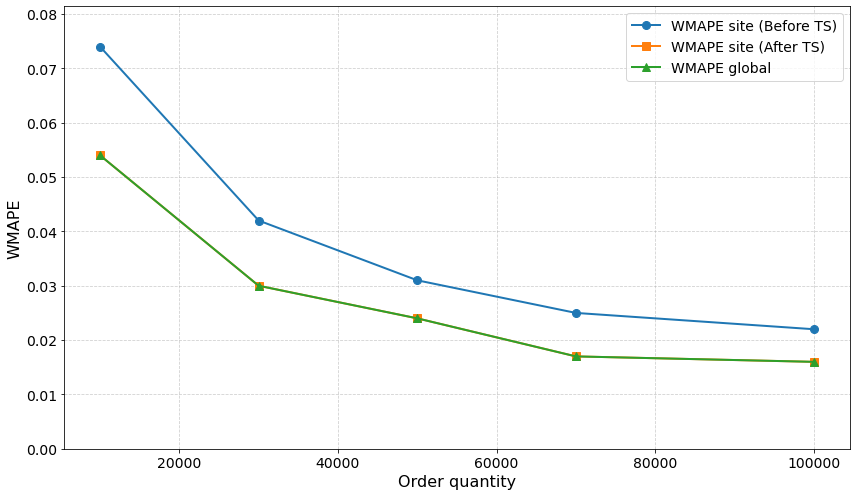}
        \caption{WMAPEs of TS}
        \label{fig:sca1}
    \end{minipage}
\end{figure}

In this test, each method is applied to five different order quantities, ranging from 10,000 to 100,000, over two days: LD12 and LD11. Figure~\ref{fig:sca} presents the optimization times of the three algorithms, while Figures~\ref{fig:sca01}, \ref{fig:sca11}, and \ref{fig:sca1} show the gap between site-level WMAPE and global WMAPE for each method. The results confirm that B\&B consistently outperforms the heuristics, delivering optimal solutions across all quantities in the shortest time. An important observation is that both errors decrease as the total order quantity increases. This finding aligns with the principles of aggregate forecasting, which suggest that forecasts are more accurate for groups of items than for individual items, as variability within a group tends to cancel out, producing more stable and reliable predictions \cite{Reid2012}. Moreover, with larger order quantities, individual fluctuations exert a smaller influence on overall error, whereas in smaller samples, changes in each order have a proportionally greater effect.

\subsubsection{Temporal fixed test}
\begin{figure}[htbp!]
    \centering
    \includegraphics[width=0.9\textwidth]{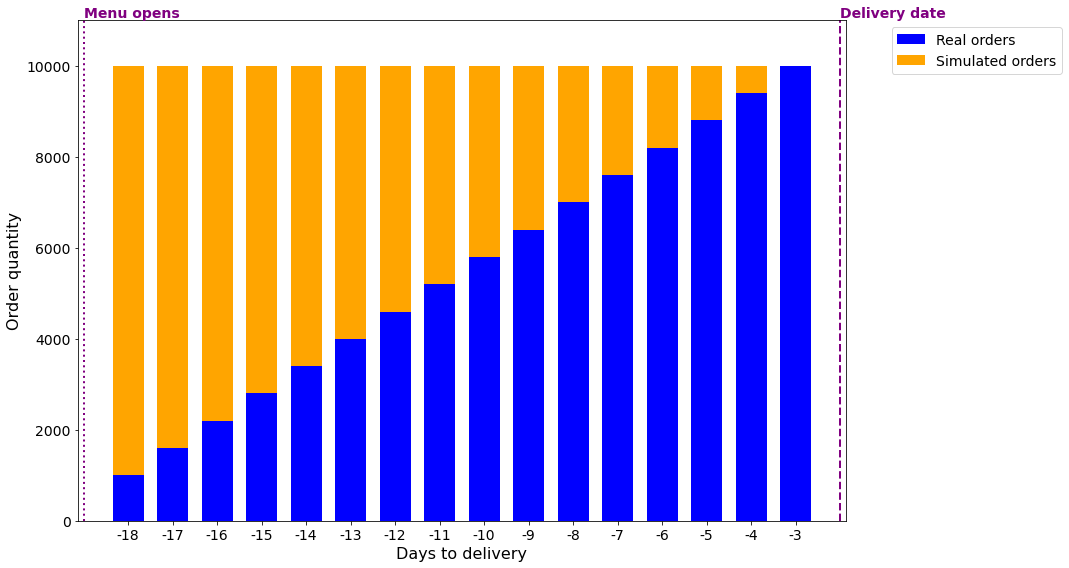}
    \caption{Composition of total orders through days}
    \label{fig:tem}
\end{figure}

In this test, the best-performing method is applied to address the temporal aspect of BAP. Specifically, B\&B is used to make allocation decisions continuously over a 15-day planning horizon under ideal conditions, assuming no changes in factory capacities or actual orders, with a fixed daily order quantity of 10,000. Figure~\ref{fig:tem} illustrates the growing proportion of real orders across days. As LD3 approaches, the order composition between consecutive days becomes increasingly similar, since each day incorporates real orders carried forward from the previous day. This explains the observed decreasing trend in the error graphs.

\begin{figure}[htbp!]
    \centering
    \begin{minipage}[b]{0.49\textwidth}
        \centering
        \includegraphics[width=\textwidth]{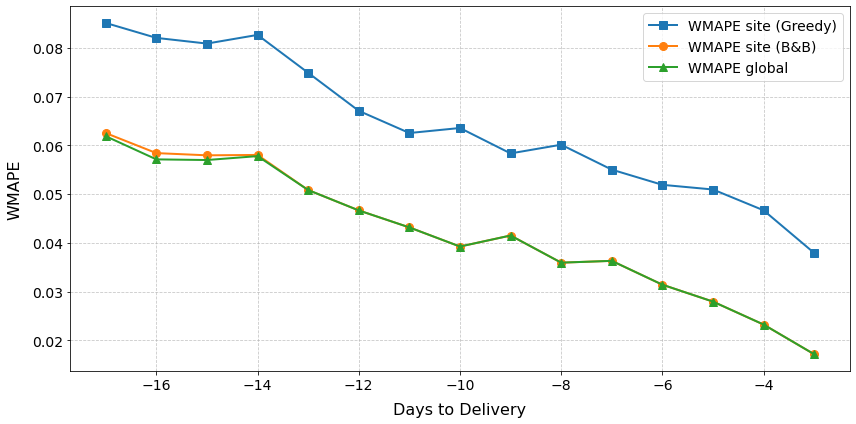}
        \caption{WMAPE site and WMAPE global through days}
        \label{fig:tem1}
    \end{minipage}
    \hfill
    \begin{minipage}[b]{0.49\textwidth}
        \centering
        \includegraphics[width=\textwidth]{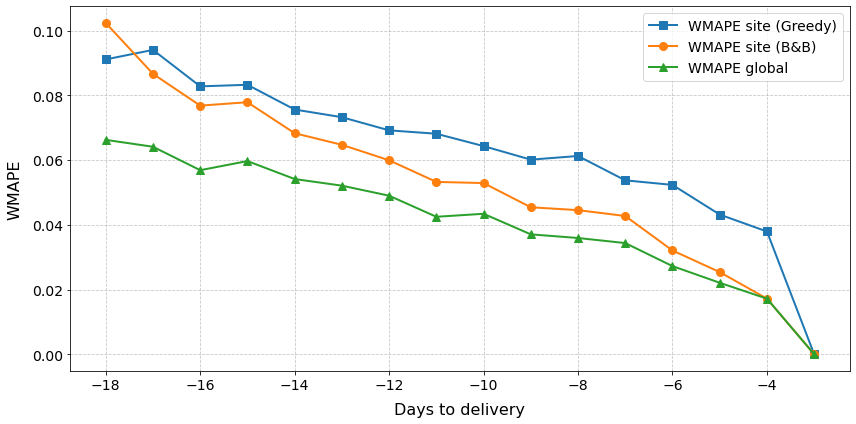}
        \caption{Each day's allocation versus LD3's allocation}
        \label{fig:tem2}
    \end{minipage}
\end{figure}

Figure~\ref{fig:tem1} shows WMAPE site when B\&B continuously optimizes allocation between two consecutive days. For instance, LD16 is optimized based on the B\&B allocation of LD17, LD15 is optimized based on LD16’s allocation, and so on. It is evident that, except for some initial days when the proportion of real orders is still low, B\&B consistently provides optimal solutions from LD14 to LD3. The WMAPE site (Greedy) line illustrates that, without B\&B, the recipe differences between two consecutive days are significantly higher.

After B\&B completes the order allocation for all days, each day’s decision is compared with the final day’s allocation, as shown in Figure~\ref{fig:tem2}. This comparison assesses whether the day-by-day optimization provides a smooth transition from soft allocations to hard allocation. The gradual decrease in WMAPE site indicates that B\&B has effectively contributed to stable production planning, minimizing abrupt changes in the allocated recipe quantity. The ultimate goal is to reduce the area under the WMAPE site curve, ideally aligning it with the area under the WMAPE global curve. Although there is a big gap between the WMAPE global and WMAPE site (B\&B) lines in the early days, this gap narrows over time. As a result, the area under the WMAPE site (B\&B) curve is smaller than that under the WMAPE site (Greedy) curve, showing B\&B’s effectiveness in achieving proxy optimization.

\subsubsection{Temporal variation test}
In the previous test, factory capacities and real orders were assumed to remain constant throughout the entire planning period. However, in practice, the following scenarios may occur:

\begin{enumerate}
    \item \textbf{Capacity change:} When F1’s capacity drops from 3,000 to 1,000 orders on LD10, there is a sudden increase in WMAPE site, as shown in Figure~\ref{fig:var11}. This spike occurs because orders must be rapidly redistributed among factories. B\&B handles this effectively, quickly adjusting allocations and achieving the optimal WMAPE site by LD9. In contrast, the greedy method exhibits a large gap between WMAPE site and global after the capacity change, highlighting its lack of adaptability. Comparing each day’s allocations with the final day’s allocation in Figure~\ref{fig:var1} shows that the WMAPE site under B\&B gradually approaches the optimal value, particularly from LD7 onward. The sharp drop in error between LD10 and LD11 reflects a significant shift in allocation decisions, with pre-LD10 allocations differing markedly from LD3. This test demonstrates B\&B’s ability to effectively respond to major changes in factory capacity.

\begin{figure}[htbp!]
    \centering
    \begin{minipage}[b]{0.49\textwidth}
        \centering
        \includegraphics[width=\textwidth]{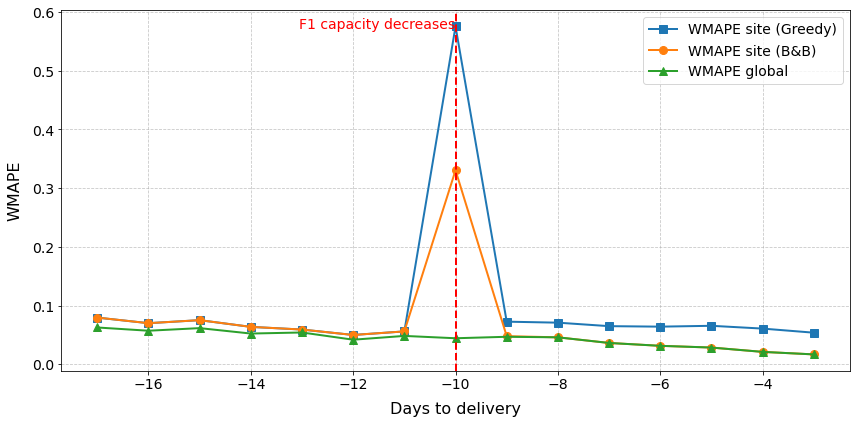}
        \caption{WMAPE site and WMAPE global through days}
        \label{fig:var11}
    \end{minipage}
    \hfill
    \begin{minipage}[b]{0.49\textwidth}
        \centering
        \includegraphics[width=\textwidth]{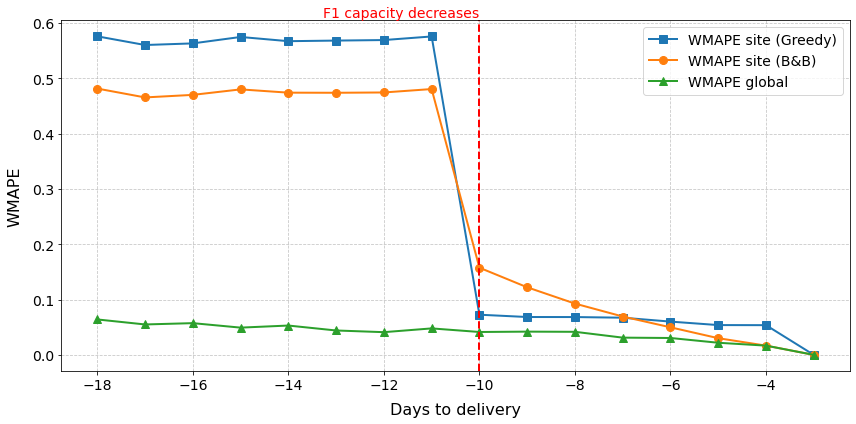}
        \caption{Each day's allocation versus LD3's allocation}
        \label{fig:var1}
    \end{minipage}
\end{figure}
    
    \item \textbf{Order changes:} This study assumes that each day, customers delete 5\% of real orders and modify recipes in 30\% of them. Figure~\ref{fig:var13} illustrates how these changes induce fluctuations in WMAPEs. Without such modifications, both WMAPE site and global exhibit a clear downward trend. When changes occur, however, their values fluctuate considerably depending on the extent of adjustments. Despite this variability, the results confirm that B\&B consistently achieves the optimal WMAPE site in all scenarios.

\begin{figure}[htbp!]
    \centering
    \includegraphics[width=0.8\textwidth]{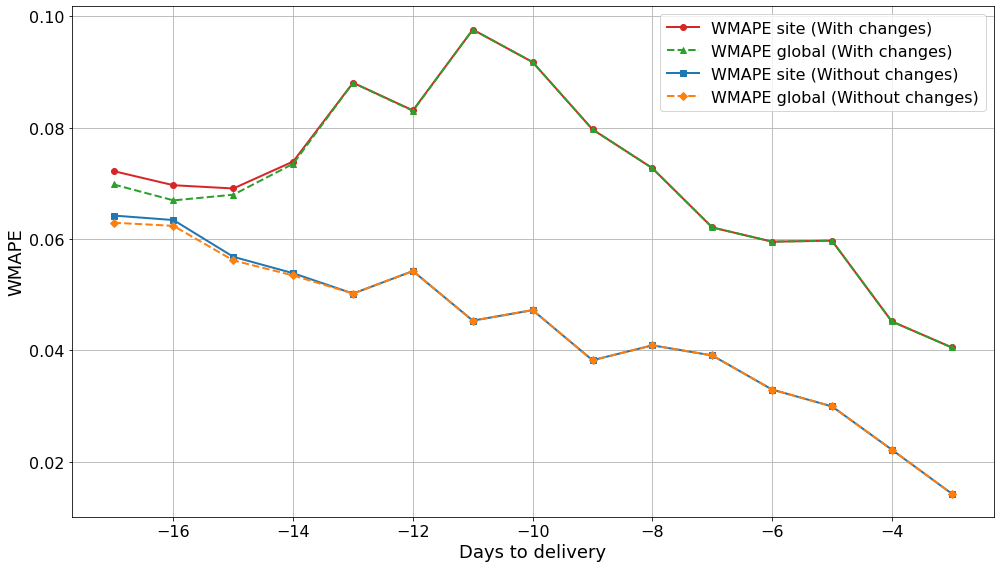}
    \caption{WMAPE site and WMAPE global with and without order changes}
    \label{fig:var13}
\end{figure}

Also in this scenario, a simple allocation strategy can be applied since the IDs of real orders remain unchanged across days. For example, if order O1 is assigned to F1 and O2 to F2 on LD14, the same assignments can be carried over to LD15. This ID-based approach ensures consistency and minimizes allocation changes for recurring orders. However, challenges emerge when recipes within an order change, potentially making it ineligible for its original factory. Deleted orders add further complications, as past real order IDs may no longer exist. Figures~\ref{fig:var23} and~\ref{fig:var3} illustrate these challenges by comparing the ID-based method with B\&B. The results show that B\&B consistently delivers optimal solutions, aligning WMAPE site with WMAPE global regardless of order changes. In contrast, even under stable condition, the ID-based method aligns with WMAPE global only during later stages, when the proportion of real orders is high. Once order changes occur, the gap between WMAPE site (ID-based) and WMAPE global widens, underscoring the limitations of simple allocation method and the advantage of more advanced techniques like CBC solver.

\begin{figure}[htbp!]
    \centering
    \begin{minipage}[b]{0.49\textwidth}
        \centering
        \includegraphics[width=\textwidth]{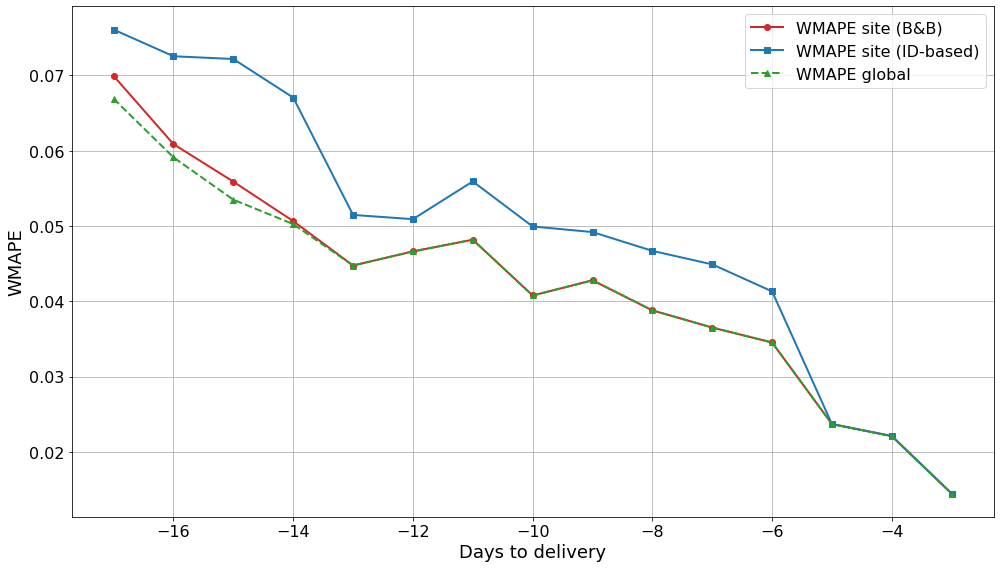}
        \caption{B\&B and ID-based allocation without order changes}
        \label{fig:var23}
    \end{minipage}
    \hfill
    \begin{minipage}[b]{0.49\textwidth}
        \centering
        \includegraphics[width=\textwidth]{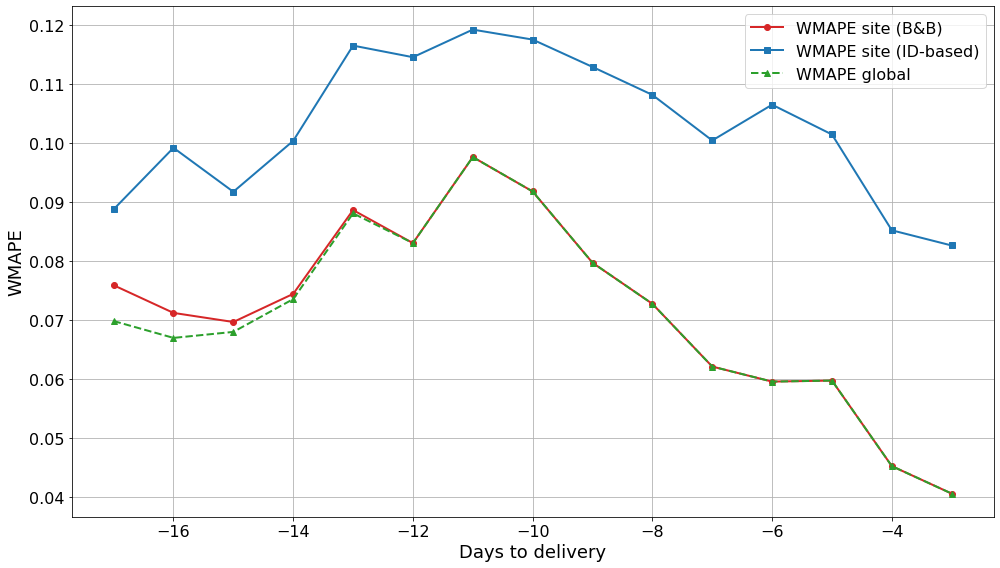}
        \caption{B\&B and ID-based allocation with order changes}
        \label{fig:var3}
    \end{minipage}
\end{figure}
\end{enumerate}

\section{Conclusion}\label{sec5}
This study introduces BAP, a novel optimization challenge in meal kit delivery operations, offering both theoretical and practical insights for its resolution. Leveraging structural similarities with BPP and CVRP, we developed a MILP formulation and systematically evaluated exact versus heuristic solution approaches. Computational experiments show that CBC consistently delivers optimal solutions with remarkable efficiency and robustness under dynamic conditions, outperforming heuristic methods such as ITPS and TS.

Three significant findings emerge from this research. First, the inverse relationship between problem scale and allocation error indicates that consolidating allocation cycles can enhance operational performance. Second, the temporal dimension of BAP imposes unique constraints not captured by static models, underscoring the importance of time-consistent allocation strategies. Third, CBC’s superior performance in dynamic settings confirms the value of advanced optimization techniques over simple rule-based methods for practical implementation.

This study is limited by its focus on consecutive-day allocations. Future work could extend BAP to multi-day horizons, integrate demand uncertainty, and explore metaheuristics such as Genetic Algorithms or Large Neighborhood Search for richer neighborhood structure. Incorporating real-time data and collaborative supplier mechanisms could further improve adaptability and long-term strategic planning, ultimately enhancing efficiency, reducing food waste, and supporting sustainability objectives in meal kit delivery system.

\section*{ACKNOWLEDGMENTS}
All numerical experiments were implemented in Python and performed on a personal computer with an Intel(R) Core(TM) i5-10210U CPU @ 1.60GHz and 8 GB of RAM. This research received financial support from the University of Southampton. The funding source had no involvement in study design; in the collection, analysis and interpretation of data; in the writing of the report; and in the decision to submit the article for publication.

\section*{DATA AVAILABILITY}
The data that support the findings of this study are openly available in GitHub at \url{https://github.com/35436506/recipe_box_allocation_optimization_in_meal_kit_delivery}.

\section*{CONFLICT OF INTEREST}
The authors declare no potential conflict of interests.

\section*{ORCID}
Thi Minh Thu Nguyen \url{https://orcid.org/0009-0000-9058-3913}\\[1ex]
Alain Zemkoho \url{https://orcid.org/0000-0003-1265-4178}

\renewcommand{\refname}{REFERENCES}

\end{document}